\documentclass[a4paper,12pt]{article}
\usepackage{amssymb}
\usepackage{amsmath}
\usepackage{graphicx}
\usepackage{epstopdf}
\usepackage{esint}

\def\Lmin{\mathcal{L}^-}

\newcommand{\eneqa}{\end{eqnarray}}
\newcommand{\begeqaet}{\begin{eqnarray*}}
\newcommand{\eneqaet}{\end{eqnarray*}}
\newcommand{\be}{\begin{equation}}
\newcommand{\ee}{\end{equation}}

\newcommand{\rn}{\rbig^n}

\newcommand{\rbig}{{\mathbb{R}}}

\newcommand{\odo}{\Omega_{d_0}}

\newcommand{\lpr}{{L^\prime}}
\newcommand{\ed}{\end{document}}

\newcommand{\dd}{\mathcal{D}}

\newcommand\ld{\mathcal{L}_D}

\newcommand\ldu{\ld[u]}

\newcommand\lr{\mathcal{L}_\rho}

\newcommand\lo{\mathcal{L}_0}

\def\beq{\begin{equation}}
\def\eeq{\end{equation}}

\def\Box{\hfill\framebox(0.25,0.25){}}

\newtheorem{thm}{Theorem}[section]

\newtheorem{prop}{Proposition}[section]

\newtheorem{lem}{Lemma}[section]

\newcommand{\begeqa}{\begin{eqnarray}}

\author{Fiorella Garcia}
\author{Boyan Sirakov}
\author{Mayra Soares}

\begin{document}


\begin{center}{\bf\Large Boundary weak Harnack estimates and regularity for elliptic PDE in divergence form}
\end{center}\smallskip

 \begin{center}
Fiorella Rend\'on, Boyan Sirakov, Mayra Soares
\end{center}\bigskip

{\small \noindent{\bf Abstract}. We obtain a global extension of  the classical weak Harnack inequality which extends and quantifies the Hopf-Oleinik boundary-point lemma, for  uniformly elliptic equations in  divergence form. Among the consequences is a boundary gradient estimate, due to Krylov and well-studied for non-divergence form equations, but completely novel in the divergence framework. Another consequence is a new more general version of the Hopf-Oleinik lemma. }
%

\section{Introduction}


We study boundary estimates and  global extensions of the  weak Harnack inequality for  PDE driven by a general
linear  uniformly elliptic second order operator in divergence form
\begin{equation}\label{defdiv}
\ldu:=-\mathrm{div}(A(x)Du +  \beta(x)u) + b(x) Du +c(x) u.
\end{equation}
 We  assume that the matrix $A$ is bounded and uniformly positive; the lower-order coefficients belong to Lebesgue spaces which make possible for weak solutions to satisfy the maximum
 principle and the Harnack inequality:
$$
\begin{array}{lc}
(H1)&A(x)\in L^\infty(\Omega),
\qquad \lambda I\le A(x)\le \Lambda I\quad \mbox{for some } 0<\lambda \le \Lambda,\\
(H2)&\hbox{$\beta,b,g\in L^q_{\mathrm{loc}}(\Omega)$ for some $q>n$, \ \ $c,f\in L^{p}_{\mathrm{loc}}(\Omega)$ for some $p>n/2$},
\end{array}
$$
for some bounded $\Omega\subset\rn$; and consider (in)equalities in the form
\begin{equation}\label{theeq}
\ldu\le(\ge,=) f+\mathrm{div}(g)\quad\mbox{ in }\;\Omega
\end{equation}
satisfied in the usual weak Sobolev sense (see  \cite[Chapter 8]{GT}) by  $u\in H^1_{\mathrm{loc}}(\Omega)$.

We  recall the De Giorgi-Moser "weak Harnack inequality" (WHI), a fundamental result in the theory of elliptic PDE. In its classical form it states that for any nonnegative supersolution of \eqref{theeq} and $B_{2R}=B_{2R}(x_0)\subset\Omega$,
\begin{equation}\label{weakha}\left(\int_{B_{R}} u^\epsilon\,dx\right)^{1/\epsilon}\le C_0\left( \inf_{B_R} u + \|f\|_{L^{p}(B_{2R})}+ \|g\|_{L^{q}(B_{2R})}\right),
\end{equation}
where  $\epsilon<n/(n-2)_+$, $C_0$ depends on $n,\lambda, \Lambda, p,q, R, \epsilon$, and the above Lebesgue norms of the coefficients $\beta, b,c$ in $B_{2R}$ (see \cite[Th. 8.18]{GT} and \cite{Tr2}).
The  WHI has a wide range of applications, the best-known being the local H\"older regularity and bounds for solutions of $\ldu= f+\mathrm{div}(g)$. As for global bounds,  essential in the study of boundary value problems, it is known that the WHI applied to $\Omega=\rn$ and  the supersolution
$u_m(x)=\min\{u(x),m\}$ if $x\in \Omega$, $u_m(x)=m$ if $x\not\in\Omega$, with $m=\inf_{\partial\Omega}u$, is sufficient to obtain global H\"older estimates if $\Omega$ has for instance the exterior cone property - see \cite[Section 8.10, Theorems 8.26, 9.27]{GT}, as well as \cite{KS2}, \cite{ARV},  for variants of this "boundary weak Harnack inequality" (bWHI).

Note this bWHI is void if $u$ vanishes on the boundary ($m=0$). One may wonder whether there is a way to quantify the positivity of the supersolution close to $\partial \Omega$ in the same way as the WHI quantifies the positivity of $u$ in the interior. Such results have appeared only recently (a review is given below), for equations in non-divergence form. Our first statement deals with this question in the previously unstudied divergence setting of \eqref{defdiv}.

We set $d(x)=\mathrm{dist}(x,\partial\Omega)$, $\odo=\{x\in\Omega\::\: d(x)<d_0\}$, and  assume that
\smallskip

\noindent{\it (H3) the boundary of $\Omega$ is $C^{1,Dini}$, the coefficients $A,\beta,g$ have Dini mean oscillation in  $\Omega_{d_0}$, and $b,c,f\in L^q(\odo)$, for some $q>n$, $d_0>0$. }
\smallskip

Below we recall and discuss the regularity notions  in (H3).
 For  $B_R=B_R(x_0)$, $x_0\in \partial\Omega$, $R\le d_0/2$, we denote $B_R^+=B_R\cap\Omega$, $B_R^0=B_R\cap\partial\Omega$, and set $k_{2R}= \|f\|_{L^q(B_{2R}^+)}+\|g\|_{L^\infty(B_{2R}^+)} + \mathcal{M}_g(B_{2R}^+)$ (here $\mathcal{M}_g$ quantifies the  Dini mean oscillation of $g$ through the function $\varrho_{m\dd,g}$  defined below). If $g\in C^\alpha$ we can take $\mathcal{M}_g$ to be the standard $\alpha$-H\"older bracket of $g$. Since all hypotheses are preserved under $C^{1,Dini}$-changes of variable, we can assume $B^0$ is flat.

\begin{thm}\label{bwhi} 1) Assume $(H1)$-$(H3)$ and  $\ldu\ge f+\mathrm{div}(g)$,  $u\ge 0$ in $B_{2R}^+$. Then for $\varepsilon>0$ depending on $n$, $\lambda$, $\Lambda$, $q$, $\varrho_{m\dd,A}$, and $C>0$ depending on $n$, $\lambda$, $\Lambda$, $q$, $\varrho_{m\dd,A}$, and  $R$,  $\|b\|_{L^q(B_{2R}^+)}$, $\|c\|_{L^q(B_{2R}^+)}$, $\|\beta\|_{L^\infty(B_{2R}^+)}$, $\mathcal{M}_\beta(B_{2R}^+)$,
 \begin{equation}\label{ineq_bwhi}
\left( \int_{B_{R}^+} \left(\frac{u}{d}\right)^\varepsilon\right)^{1/\varepsilon}\le C\left(\inf_{B_R^+} \frac{u}{d} + k_{2R}\right).
\end{equation}

2) Assume $\ld^{(1)}, \ld^{(2)}$ are operators in the form \eqref{defdiv} under  $(H1)-(H3)$ and  $\ld^{(1)}[u]\le f^{(1)}+\mathrm{div}(g^{(1)})$, $\ld^{(2)}[u]\ge f^{(2)}+\mathrm{div}(g^{(2)})$, $u\ge 0$ in $B_{2R}^+$, $u=0$ on $B_{2R}^0$. There exists  $C>0$ depending on $n$, $\lambda$, $\Lambda$, $q$, $\varrho_{m\dd,A^{(i)}}$, and $R$,  $\|b^{(i)}\|_{L^q(B_{2R}^+)}$, $\|c^{(i)}\|_{L^q(B_{2R}^+)}$, $\|\beta^{(i)}\|_{L^\infty(B_{2R}^+)}$, $\mathcal{M}_{\beta_i}(B_{2R}^+)$, $i=1,2$, such that
 \begin{equation}\label{ineq_bhi}
\sup_{B_R^+} \frac{u}{d}\le C\left(\inf_{B_R^+} \frac{u}{d} + k_{2R}^{(1)}+ k_{2R}^{(2)}\right).
\end{equation}

3) All $B^+$ (resp. $B^0$) can be replaced by $\Omega$ (resp. $\partial\Omega$) in 1) and 2), with $C$ depending also on $p$ and $\Omega$.
\end{thm}

\noindent{\it Remark.} The functions $\varrho$ are defined in the next section. For a coefficient $\xi\in C^\alpha$, $\alpha\in(0,1)$, the dependence of $\varepsilon, C$ in $\varrho_{m\dd,\xi}$ (resp. in $\mathcal{M}_\xi$) reduces to dependence in $\alpha$ and the $C^\alpha$-norm of the coefficient.\medskip

The classical and fundamental Zaremba-Hopf-Oleinik lemma, or ``boundary point principle" (BPP), states that a nontrivial nonnegative supersolution of $\mathcal{L}[u]\ge0$ is indeed ``strictly positive" close to $\partial \Omega$, in the sense that $\inf_\Omega(u/d)>0$, for a sufficiently smooth $\partial\Omega$.
A lot of work has been dedicated to  getting optimal conditions for the validity of BPP, in terms of the regularity or the geometry of the domain, or of the nature of the coefficients of the elliptic operator. We refer to  \cite{AM},  \cite{AN2}, \cite{AN3}, \cite{CLN}, \cite{FG},  \cite{K1}, \cite{LXZ}, \cite{N1},   \cite{S1},  \cite{S2} for such conditions, as well as historical reviews and more references. The survey \cite{AN3} is very complete and up-to-date.

Note that \eqref{ineq_bwhi} with $f=g=0$ ($k_{2R}=0$) implies the BPP, and quantifies it in the following sense: if $\ldu\ge0$ in $B_2^+$, $u\ge  c_0d$ in a subset $\omega\subset B_1^+$ of positive measure, then $u\ge\kappa c_0d$ in the whole $B_1^+$, for some $\kappa>0$ depending only on $|\omega|$ and the data -- this can be thought of as a boundary variant of a ``growth lemma".  To our knowledge, there are no previous results which quantify the BPP in such a way for any type of divergence form equations which cannot be related to non-divergence ones. In particular, \eqref{ineq_bwhi} is new for inequalities such as $-\mathrm{div}(A(x)Du)\ge0$ or $-\Delta u \ge \mathrm{div}(g)$, with $A,g\in C^\alpha$.

  Furthermore, even the BPP itself implied by Theorem~\ref{bwhi} with $f=g=0$ appears to be new for non-Dini continuous leading coefficients. The best available hypothesis under which the BPP was proved for $-\mathrm{div}(A(x)Du)\ge0$, was Dini continuity of $A$, see \cite{AN3}.
  Previous works on the BPP also make  hypotheses on the distributional sign (or absence) of div$(\beta)-c$.

The importance of such a quantification of the BPP was recognized only recently, but already a number of applications have appeared. The uniform up-to-the boundary  inequality \eqref{ineq_bwhi} for non-divergence form equations was proved in \cite{Sir1}.  In the non-divergence framework we also refer to \cite[Lemma 1.6]{CLN},  \cite{BM}, \cite{LXZ}, for estimates like \eqref{ineq_bwhi} in which the left-hand side contains an integral on a interior subset, and the constant $C$ degenerates if this subset approaches the boundary.  The best constant $\varepsilon$ for which \eqref{ineq_bwhi} holds was specified in \cite{Sir3}, for operators which are both in divergence and non-divergence form (the optimal value of $\varepsilon$ is an open problem for more general operators). The results in \cite{Sir1}, \cite{Sir3} have been instrumental in a new method for a priori bounds for positive solutions of nonlinear elliptic equations -- see \cite{Sir2} and the references there. Another application has just appeared in \cite{KS2022}. We will use Theorem \ref{bwhi} in the boundary regularity Theorem \ref{kryl} below, as well as in  the forthcoming works \cite{GS} on  solvability  of equations having quadratic dependence in the gradient and \cite{SS2} on the Landis conjecture and elliptic estimates with optimized constants.

As often happens in elliptic theory, the statement of Theorem \ref{bwhi} is similar to that of the non-divergence case
\cite[Theorem 1.2]{Sir1};
however, the main point of the proof (the boundary growth lemma) requires a different approach. Here we use the classical idea of \cite{FG} to compare $u$ with a solution of a ``frozen coefficients" equation in a sufficiently small annulus which touches the boundary; however, we combine this comparison with direct use of elliptic estimates, in particular the Stampacchia maximum principle and the global $C^1$-estimates from \cite{DEK}, thus avoiding the  use of Green functions which has been frequent in proofs of the BPP in the  divergence framework. 
\medskip

Our second main result is an application of Theorem \ref{bwhi} to  boundary regularity theory. It concerns the following classical property: {\it given two elliptic operators such that the solutions of the Dirichlet problem in $\Omega$ for each of them have uniformly continuous gradient in $\overline{\Omega}$; and a function which is only a subsolution and a supersolution of two} different {\it equations involving these operators, then this function may not even be differentiable in $\Omega$ but still has a uniformly continuous gradient at $\partial \Omega$}. This is a fundamental result in the non-divergence theory, which goes back to Krylov and his proof of solvability and regularity of the Dirichlet problem \cite{K1}. Krylov's property  has been studied, extended and used over the years by many authors, see \cite{S0}, \cite{LU2}, \cite{LinLi}, \cite{SilS}, \cite{LZ}, and \cite{BGMW} for very general results and a large discussion, as well as the references in these works. However, this fact has never been proven for pure divergence-form equations, even in the simplest cases.

 We show that Krylov's property is valid in the divergence framework.

 \begin{thm}\label{kryl} Assume $\partial \Omega$ is in $C^{1,\overline \alpha}$, $\ld^{(1)}, \ld^{(2)}$ are operators in the form \eqref{defdiv} under  $(H1)$, whose coefficients $A^{(i)},\beta^{(i)},g^{(i)}\in C^{\overline \alpha}(B_{1}^+)$, $b^{(i)},c^{(i)},f^{(i)}\in L^q(B_{1}^+)$, for some $\overline \alpha>0$, $q>n$, $i=1,2$. Assume $u\in H^1(B_{1}^+)$ is such that  $$\ld^{(1)}[u]\le f^{(1)}+\mathrm{div}(g^{(1)}),\quad \ld^{(2)}[u]\ge f^{(2)}+\mathrm{div}(g^{(2)})\mbox{ in }B_{1}^+,\quad  u|_{B_{1}^0}\in C^{1,\overline \alpha}(B_{1}^0).$$
Then there exists $G \in C^{\alpha}(B_{1/2}^0,\rn)$ (the``gradient" of $u$ on $B_{1/2}^0$),  such that
\begin{equation} \label{e:targetC1a0}
\|G\|_{C^{\alpha}(B_{1/2}^0)} \leq CW,
\end{equation}
and for every $x  \in B_{1/2}^+$ and every $x_0\in B_{1/2}^0$ we have
\begin{equation} \label{e:targetC1a}  |u(x) - u(x_0)-  G(x_0) \cdot (x-x_0)  | \leq C W |x-x_0|^{1+\alpha},\quad\mbox{
where}
\end{equation}
$$
W:=\|u\|_{L^\infty(B_{1}^+)}+\|u\|_{C^{1,\alpha}(B_{1}^0)}+L, \qquad L=\sum_i (
\|f^{(i)}\|_{L^q(B_{1}^+)}+\|g^{(i)}\|_{C^{\overline \alpha}(B_{1}^+)}).
$$
Here $\alpha,C>0$  depend on $n$, $\lambda$, $\Lambda$, $q$,  $\overline \alpha$, $\|A^{(i)}\|_{C^{\overline \alpha}(B_{1}^+)}$; $C$ also depends on  the H\"older, resp. Lebesgue, norms of the lower-order coefficients and  $\partial\Omega$.
\end{thm}

In Theorem \ref{kryl} we strengthened the regularity assumptions on the coefficients to the most important and often encountered H\"older continuity. This permits to us to ease technicalities and present the result as a consequence from Theorem~\ref{bwhi} and the method developed in \cite{SilS} for the non-divergence case.

In  the next section we give some more comments on our hypotheses and framework. The last section is devoted to the proofs of the theorems.

\section{Further comments}

The distinction ``divergence" vs. ``non-divergence" is particularly relevant and delicate with regard to the BPP and its ramifications.
For nondivergence type inequalities, say tr$(A(x)D^2u)\le0$, the BPP is true for any $A(x)\in L^\infty(\Omega)$. On the other hand, for inequalities in divergence form, say div$(A(x)Du)\le0$, the BPP may fail even for $A(x)\in C(\overline{\Omega})$ (see  \cite{N1}, \cite{AN3}, for counterexamples and more references). However,  the BPP is true for that inequality if $A$ is Dini continuous (see \cite{AN2},\cite{AN3}), and as we now know by Theorem~\ref{bwhi},  even if $A$ has Dini mean oscillation. 
Furthermore, it is rather remarkable that the standard  boundary Harnack inequality (in which two positive solutions are compared close to the boundary, as opposed to one solution and the distance function) is valid for div$(A(x)Du)=0$ with $A\in L^\infty$ (see \cite{CFMS}), but fails for div$(A(x)Du)=f$, $f\not=0$; however for the latter it is true if $A$ is only continuous,  as was recently shown in \cite{RT}. Another example of how delicate the role of the regularity assumptions on the coefficients may be are the recent deep works on ``propagation of smallness" (see \cite{LM})  for solutions of div$(A(x)Du)=0$, which are valid for a symmetric Lipschitz $A$, but fail for $A\in C^\alpha$, $\alpha<1$. In a certain sense, Theorem \ref{bwhi} above is a "propagation of smallness" of  $u/d$ from the boundary to the whole of the domain.\medskip

Next we comment on the regularity and integrability assumptions we make on the coefficients of the elliptic operators. We crucially use that the standard Dirichlet problem associated to the operator has global $C^1$-estimates. The Dini mean oscillation assumption on the leading coefficients is currently the most general available hypothesis under which a $C^1$-estimate up to the boundary is known, while mere continuity is not sufficient for such an estimate; we believe our method is sufficiently versatile to adapt to other situations, if global $C^1$-estimates are proved in the future under even more general assumptions.

We have assumed that the lower-order coefficients in $\ld$ belong to $L^q$ with $q>n$, which is certainly the optimal Lebesgue integrability for Theorem \ref{bwhi} (and even for the BPP, which is known to fail for instance for $b\in L^n$, see \cite[Example 4.1]{S1}). On the other hand, the BPP is known under finer restrictions on the lower-order coefficients, such as intermediate spaces between $L^q$ for $q>n$, and $L^n$,  see \cite{AN2}, \cite{AN3}, and the references there.
Theorem \ref{bwhi} should be true under such assumptions too; however, since it is new even for operators without lower order coefficients, and to avoid technical complications, we do not study such extensions here. Our assumption permits to us to directly quote the $C^1$-estimate in \cite[Theorem 1.3]{DEK} and concentrate on its use.

Similarly, while  Theorem \ref{kryl} is proved in large generality (and is new in the simplest cases such as equations without lower-order coefficients and with zero right-hand side),  we expect and conjecture that it is true  for even more general coefficients and operators. It should be possible to replace the H\"older by Dini mean continuity in the assumptions on the leading coefficients; however this would render the rescaling argument which is in the core of the proof considerably more delicate. Furthermore, the result should be true for quasi-linear operators whose associated Dirichlet problem has $C^1$ estimates, such as operators considered in \cite{LU}. For instance, we expect Theorem \ref{kryl} to be valid for operators with quadratic growth in the gradient as in \cite{BGMW}, replacing the Pucci operators there by divergence form operators with coefficients which have  Dini mean oscillation.

\section{Proofs}

\subsection{Preliminaries}
We start by recalling the $C^1$ estimate from \cite{DEK}. Following that paper, a function $\varrho:[0,1]\to[0,\infty]$ is a Dini function (we write $\varrho\in\dd$) if $\varrho(0)=0$, $0<c_1\rho(t) \le \rho(s)\le c_2\rho(t)$ for $0<t/2\le s\le t$
 and $\int_0 (\varrho(s)/s)\,ds$ converges. We note it is possible to assume without restricting the generality that $\varrho(s)$ is non-decreasing and continuously differentiable for $s>0$, and $\varrho(s)/s$ is non-increasing, see \cite{AN3} (so we can take $c_1=1/2$, $c_2=1$).

A function $h$ is Dini continuous on $\Omega$ (we write $h\in C^{\dd}(\Omega)$) if
$$
\varrho_{D,h}(r) = \sup_{x\in\Omega}\sup_{y^\prime, y^{\prime\prime}\in B_r^+(x)} |h(y^\prime)-h(y^{\prime\prime})|\in \dd.
$$
We say that $\Omega$ is in $C^{1,\dd}$ if each point on $\partial\Omega$ has a neighborhood in which $\partial\Omega$ is the graph of a continuously differentiable function whose derivatives are in $C^{\dd}$.

A function $h$ has Dini mean oscillation on $\Omega$ (we write $h\in C^{m\dd}(\Omega)$) if
$$
\varrho_{m\dd,h}(r) = \sup_{x\in\Omega}\fint_{B_r^+(x)} \left|h(y) -\fint_{B_r^+(x)}h(y)dy\right|\,dy \in \dd.
$$
Here as usual $\fint_G= \frac{1}{|G|}\int_G$. Note that $\varrho_{m\dd,h}(r)\le \varrho_{D,h}(r)$, so Dini mean oscillation is a weaker hypothesis than Dini continuity. A standard example of non-Dini continuous function which has Dini mean oscillation is $h(x) = |\log|x||^{-\gamma}$, $\gamma\in (0,1]$ (for  enlightening examples on the difference between Dini and mean Dini conditions, see \cite{DK}, \cite{MM}). For $k>0$, $t\in(0,1)$, under the rescaling \begin{equation}\label{resc}\tilde{h}(y) = kh(ty),\;\mbox{ we have } \; \varrho_{\dd,\tilde{h}}(r)=k\varrho_{\dd,{h}}(tr), \quad \varrho_{m\dd,\tilde{h}}(r)=k\varrho_{m\dd,{h}}(tr).
\end{equation}
 It is also true that in a $C^1$--domain $\Omega$, if $h\in C^{m\dd}(\Omega)$ then $h$ is uniformly continuous in $\Omega$, with a modulus of continuity $\omega_h(r)$ dominated by $\int_0^r (\varrho_{m\dd,h}(s)/s)\,ds$ (see \cite[Lemma A.1]{HK}).

\begin{thm}\label{dek} (Dong-Escauriaza-Kim, \cite{DEK}) Assume $(H1)$ and $(H3)$ hold in~$\Omega$, $\partial\Omega\in C^{1,\dd}$, $\mathrm{diam}(\Omega)\le1$. If $u\in H^1_0(\Omega)$ solves $\ldu= f+\mathrm{div}(g)$ in~$\Omega$, then $u\in C^1(\overline{\Omega})$. In addition,
\begin{equation}\label{c1}\|u\|_{C^1(\Omega)}\le C (\|u\|_{L^2(\Omega)}+ \|f\|_{L^q(\Omega)}+  \|g\|_{L^\infty(\Omega)}+\mathcal{M}_g(\Omega)),
\end{equation}
where $\mathcal{M}_g(\Omega)$ is a quantity which describes the Dini mean  oscillation of $g$ and is defined through the values of $\varrho_{m\dd,g}(r)$ (in particular of $\int^r_0 (\varrho_{m\dd,g}(s)/s)$, $r\in(0,1]$). The constant $C$ is bounded above in terms of $n,\lambda, \Lambda, q$, the $C^{1,\dd}$-norm of $\partial\Omega$,  upper bounds on the $L^q$-norms of $b,c$, the $L^\infty$-norm of $\beta$, $\mathcal{M}_A(\Omega)$, $\mathcal{M}_\beta(\Omega)$.

 Furthermore, there exists a modulus of continuity $\sigma$ determined by $n,\lambda, \Lambda$, $q$,  the $L^q$-norms of $b,c,f$,  the functions $\varrho_{mD}$, $\varrho_{\dd}$ corresponding to $A,\beta,g$ and $\partial\Omega$, and by $\|u\|_{L^2(\Omega)}$,  such that
$$
|Du(x)-Du(y)|\le \sigma(|x-y|).
$$

\end{thm}
This statement can be inferred from  \cite[Theorem 1.3]{DEK} and its proof. See in particular inequalities (2.31) and (2.36) in \cite{DEK}. Note the term $\|Du\|_{L^1(\Omega)}$ which appears in (2.31) is bounded by the right-hand side of \eqref{c1} by the standard Sobolev bounds for weak solutions (see for instance  \cite[Theorem~3.2]{St})
\begin{equation}\label{stan}\|u\|_{H^1(\Omega)}\le C (\|u\|_{L^2(\Omega)}+ \|f\|_{L^q(\Omega)}+  \|g\|_{L^q(\Omega)}).
\end{equation}

\noindent {\it Remark 3.1}.  We have that $\mathcal{M}_{tg}(\Omega)=t\mathcal{M}_{g}(\Omega)$ for $t>0$ and for each $\epsilon>0$ there is $\delta>0$ such that  $\mathcal{M}_g(\Omega)<\epsilon$ if $\int^1_0 (\varrho_{m\dd,g}(s)/s)\,ds<\delta$.

\noindent {\it Remark 3.2}. If $g\in C^\alpha$ we can take $\mathcal{M}_g$ to be the usual $\alpha$-H\"older bracket of~$g$.

\subsection{Proof of Theorem \ref{bwhi}}

We observe that all hypotheses on the operator $\mathcal{L}_D$ are preserved under a $C^{1,\mathcal{D}}$-regular change of variables (that Dini mean oscillation is preserved follows from \cite[Lemma 2.1]{DEK}). So from now on we will assume that the boundary of $\Omega$ is locally flat, included in $\{x:x_n=0\}$.
In the following $Q_\rho=Q_\rho(\rho e)$ denotes the cube with center $\rho e$ and side $\rho$, where $e= (0,\ldots,0,1/2)$. To avoid confusion, the reader's attention is brought to the fact that $Q_\rho$ is not centered at the origin but has its bottom on $\{x_n=0\}$.

 We first establish the following growth lemma. We assume that in $Q_2$ all coefficients have the regularity given in $(H3)$, and set
$$
W=\|f\|_{L^q(Q_{2})}+\|g\|_{L^\infty(Q_2)}+ \mathcal{M}_g(Q_{2}).
$$

\begin{lem}\label{deggrowth} Given $\nu>0$, there exist $a,k\in(0,1)$ depending on $n$,  $\lambda$, $\Lambda$,  $q$, $\varrho_{mD,A}$, $\nu$,   such that if $u\in H^1(Q_2)$ is a weak solution of
$$\ldu\ge  f+\mathrm{div}(g), \quad u\ge0\mbox{ in }Q_2, $$
$$
\|b\|_{L^q(Q_2)}\le 1,\quad \|c\|_{L^q(Q_2)}\le 1, \quad \|\beta\|_{L^\infty(Q_2)}+ \mathcal{M}_\beta(Q_{2})\le 1, \quad W \le a,
$$
and we have
\begin{equation}\label{condigrowth}
|\{u>x_n\}\cap Q_1|\ge \nu ,
\end{equation}
then $u> kx_n$ in $Q_1$.
\end{lem}

\noindent{\it Proof.} For all $\rho\in(0,1]$  set $x_\rho=(0,\ldots,0,\rho)$,  $A_\rho=B_\rho(x_\rho)\setminus B_{\rho/2}(x_{\rho})$. For  $y\in A_1$, $v\in H^1(A_1)$ introduce the operator
$$
\lr[v]:=-\mathrm{div}_y(A(\rho y)Dv(y) +  \beta_\rho(y)v(y)) + b_\rho(y) D_yv(y) +c_\rho(y) v(y),
$$
$$
\mbox{where}\qquad\beta_\rho(y) = \rho \beta(\rho y), \quad b_\rho(y)= \rho b(\rho y), \quad c_\rho(y)= \rho^2 c(\rho y).
$$
Note $\ldu\le  f+\mathrm{div}(g)$ in $A_\rho$ is equivalent to $$\lr[u_\rho]\le f_\rho(y) +\mathrm{div}_y(g_\rho(y))\quad\mbox{ in }\; A_1,$$ if we set $$y=x/\rho,\quad u_\rho(y) = u(x),\quad f_\rho(y) = \rho^2 f(\rho y),\quad g_\rho(y)= \rho g(\rho y).$$ Then
\begin{equation}\label{scalq}
\|\beta_\rho\|_{L^q(A_1)}= \rho^{1-n/q}\|\beta\|_{L^q(A_\rho)}\le \rho^{1-n/q}\|\beta\|_{L^q(Q_2)}\le 2^n\rho^{1-n/q},
\end{equation}
 and similarly for the other coefficients with subscript $\rho$. So all these coefficients have $L^q$ norms in $A_1$ that tend to zero as $\rho\to0$, since $1-n/q>0$ for $n<q\le\infty$. Also by the definition of $\varrho_{mD}$, \eqref{resc} and Remark 3.1 we have
 \begin{equation}\label{scalq1}\mathcal{M}_{g_\rho}(A_1)\le \rho \mathcal{M}_{g}(Q_2)\le\rho\quad\mbox{ and }\quad\mathcal{M}_{\beta_\rho}(A_1)\le \rho \mathcal{M}_{\beta}(Q_2)\le \rho.
 \end{equation}

Fix a smooth function $\psi$ such that $\psi=1$ in $B_{1/2}$,  $\psi=0$ outside $B_1$, $\|\psi\|_{C^1}= C(n)$. Set
$$
\bar f_\rho = f_\rho-b.D\psi - c\psi, \qquad \bar g_\rho = g_\rho + AD\psi + \beta\psi.
$$
It follows from classical solvability results (see Theorems 3.1 and 3.3 of \cite{St}) that for some $\rho_0\in(0,1/2)$ depending only on $n, \lambda, \Lambda, q$, and for all $\rho\in(0,\rho_0)$ there is a unique function $$w_\rho\in H^1_0(A_1)\quad\mbox{ such that }\quad
\lr[w_\rho]= \bar f_\rho +\mathrm{div}(\bar g_\rho)\quad\mbox{in }\; A_1.$$ Hence $v_\rho=w_\rho+\psi\in H^1(A_1)$ solves
\begin{equation}\label{defvrho}
\left\{
\begin{array}{rclcc}
\lr[v_\rho] &=& f_\rho +\mathrm{div}(g_\rho)&\mbox{in}& A_1\\
v_\rho&=&1&\mbox{on}& \partial B_{1/2}\\
v_\rho&=&0&\mbox{on}& \partial B_{1}.
\end{array}
\right.
\end{equation}

We set $\lo[\cdot]=-\mathrm{div}(A(0)D\cdot)$ (note $\lo$ has constant coefficients and is also in non-divergence form) and let $v_0$ be the solution of
\begin{equation}\label{defvo}
\left\{
\begin{array}{rclcc}
\lo[v_0] &=& 0&\mbox{in}& A_1\\
v_0&=&1&\mbox{on}& \partial B_{1/2}\\
v_0&=&0&\mbox{on}& \partial B_{1}.
\end{array}
\right.
\end{equation}
By the maximum principle $0< v_0< 1$ in $A_1$. Since $\lo[1]=0$, if we extend $v_0=1$ in $B_{1/2}$ we obtain a supersolution, $\lo[v_0]\ge0$ in $B_1$. By theorem 4.1.2 in \cite{Sir1} (or Theorem 1.2 in \cite{Sir3}) we have
\begin{equation}\label{eas}
v_0(y)\ge c_0\, \mathrm{dist}(y,\partial B_1)
\end{equation}
for all $y\in A_1$, and some $c_0>0$ depending only on $n,\lambda,\Lambda$. Also, by standard elliptic estimates for equations with constant coefficients
\begin{equation}\label{C2}
\|v_0\|_{C^2(A_1)}\le C_0= C_0(n,\lambda,\Lambda).
\end{equation}

Set $z_\rho=v_\rho-v_0$. We have
\begin{eqnarray*}
\lr[z_\rho]&=& f_\rho +\mathrm{div}(g_\rho) + (\lo-\lr)[v_0]\\
&=& \widetilde{f}_\rho + \mathrm{div}(\widetilde{g_\rho})
\end{eqnarray*}
in $A_1$, where
$$
\widetilde{f}_\rho = f_\rho - b_\rho Dv_0 - c_\rho v_0
$$
$$
\widetilde{g_\rho} = (A(\rho y) - A(0)) Dv_0 + \beta_\rho v_0 + g_\rho.
$$
Clearly by \eqref{scalq} and \eqref{C2}
$$
\|\widetilde{f}_\rho\|_{L^q(A_1)} \le C\rho^\alpha, \qquad \|\widetilde{g}_\rho\|_{L^q(A_1)}  \le C(\omega_A(\rho)+\rho^\alpha)
$$
(here $\alpha=1-n/q$ and $\omega_A$ is the uniform modulus of continuity of $A$ given by \cite[Lemma A.1]{HK}). Further, by writing the last equation in the form
\begin{eqnarray*}
\widehat{\lr}[z_\rho] &=& \widetilde{f}_\rho -c_\rho(y)z_\rho(y)+\mathrm{div}(\widetilde{g}_\rho - \beta_\rho(y)z_\rho(y)) \\
&= :& \widehat{f}_\rho + \mathrm{div}(\widehat{g_\rho})
\end{eqnarray*}
where
$$
\widehat{\lr}[z_\rho]=-\mathrm{div}(A(\rho y)Dz_\rho(y) ) + b_\rho(y){\cdot} Dz_\rho(y)
$$
and by applying Stampacchia generalized maximum principle  (\cite[Theorem 8.16]{GT}), since $z_\rho=0$ on $\partial A_1$ we get
\begin{eqnarray*}
\|z_\rho\|_{L^\infty(A_1)} &\le & C( \|\widehat{f}_\rho\|_{L^q(A_1)} +  \|\widehat{g}_\rho\|_{L^q(A_1)}) \\ &\le & C_1(\|\widetilde{f}_\rho\|_{L^q(A_1)} +  \|\widetilde{g}_\rho\|_{L^q(A_1)}) + C_2\rho^\alpha\|z_\rho\|_{L^\infty(A_1)}
\end{eqnarray*}
so setting $\rho_1>0$ such that $C_2\rho_1^\alpha=1/2$ we have for $\rho<\rho_2=\min\{\rho_0,\rho_1\}$
$$
\|z_\rho\|_{L^\infty(A_1)}\le C(\|\widetilde{f}_\rho\|_{L^q(A_1)} +  \|\widetilde{g}_\rho\|_{L^q(A_1)})\le C(\omega_A(\rho)+\rho^\alpha).
$$
The first of these inequalities is independent of the form of $\widetilde{f}_\rho, \widetilde{g}_\rho$, and means that the generalized maximum principle holds for the equation $\lr[z_\rho]= \widetilde{f}_\rho + \mathrm{div}(\widetilde{g_\rho})$ in $A_1$, and the comparison principle is valid for $\lr$ in $A_1$. By standard Sobolev estimates for weak solutions $z_\rho\in H^1_0(A_1)$ of this equation
$$
\|z_\rho\|_{H^1(A_1)}\le C(\|z_\rho\|_{L^2(A_1)}+ \|\widetilde{f}_\rho\|_{L^q(A_1)} +  \|\widetilde{g}_\rho\|_{L^q(A_1)}) \le C(\omega_A(\rho)+\rho^\alpha).
$$
We now apply Theorem \ref{dek} to the equation $\lr[z_\rho]= \widetilde{f}_\rho + \mathrm{div}(\widetilde{g_\rho})$. Note that $\|\widetilde{g_\rho}\|_{L^\infty(A_1)}+ \mathcal{M}_{\widetilde{g_\rho}}(A_1)\to 0$ as $\rho\to0$ by the homogeneity property in  \eqref{resc}, Remark 3.1, \eqref{scalq}, \eqref{scalq1}. Therefore, Theorem \ref{dek} implies that there exists  $\sigma_\rho$ with $\sigma_\rho\to 0$ as $\rho\to0$, such that
$$
\|z_\rho\|_{C^1(A_1)}\le \sigma_\rho.
$$
By the mean value theorem, for each $y=(0,\ldots,0,y_n)\in A_1$
$$
\left|\frac{z_\rho(y)}{y_n}\right|=\left|\frac{z_\rho(y)-z_\rho(0)}{y_n}\right| \le \sigma_\rho.
$$

Finally, we get by \eqref{eas}
$$
\frac{v_\rho(y)}{y_n}=\frac{v_0(y)}{y_n}+\frac{z_\rho(y)}{y_n} \ge  c_0-\sigma_\rho\ge c_0/2,\qquad y=(0,\ldots,0,y_n)\in A_1.
$$
provided $\rho\in(0,\rho_3)$, for $\rho_3>0$ fixed so that $\sigma_\rho\le c_0/2$ for $\rho\in(0,\rho_3)$.

By \eqref{condigrowth} there is $\rho_4>0$ depending on $\nu$  such that
$$
|\{u>\rho_4\}\cap Q_1\cap \{x_n\ge \rho_4\}|\ge \nu/2.
$$
Set $\bar\rho=\min_{0\le i\le 4} \rho_i$. Then by the interior weak Harnack inequality applied in $\{x_n\ge \bar\rho/4\}\cap Q_2$  to the inequality satisfied by $u$  we get
$$
u\ge c\rho_4\nu- C^\prime W\ge  c^\prime - C^\prime a
$$
in $\{x_n\ge \bar\rho/2\}\cap Q_1$ and hence by choosing $a$ sufficiently small $u\ge c_{\bar\rho}>0$ in that set. This implies that $u_{\bar\rho}\ge c_{\bar\rho} v _{\bar\rho}$ on $\partial A_1$, and since $\lr[u_{\bar\rho}]\ge \lr[v_{\bar\rho}]$ in $A_1$,  by the comparison principle, $v_{\bar\rho}\ge c_0/2$ in $A_1$ and $x_n=\bar\rho y_n$ we get
$$
\frac{u(x)}{x_n}\ge \min\{(c_0c_{\bar\rho}/(2\bar \rho),c_{\bar\rho}\},
$$
provided $x=(0,\ldots,0,x_n)\in Q_1$. We can shift the origin along $\{x_n=0\}$ so Lemma \ref{deggrowth} is proved.
\medskip

\noindent{\it Proof of Theorem \ref{bwhi}}. It is enough to prove
 \eqref{ineq_bwhi} for $R=1$, $x_0=0$ (the general case follows by scaling and translation). We can also assume $d_0=2$. We set
$$
r_0 = (4+\|b\|_{L^q(B_{2}^+)}^{\frac{1}{1-n/q}} + \|c\|_{L^q(B_{2}^+)}^{\frac{1}{2-n/q}}+\|\beta\|_{L^\infty(B_{2}^+)} +\mathcal{M}_\beta(B_{2}^+))^{-1}.$$

By the same change of variables as in the beginning of the proof of Lemma~\ref{deggrowth}, precisely for $\rho=r_0$, replacing $\ld$ by $\mathcal{L}_{r_0}$ and $u$ by $u_{r_0}$ we obtain an equation in a set containing $B_4^+$, whose lower order coefficients $\beta, b,c $ have the bounds required in Lemma \ref{deggrowth}. Note that if Theorem \ref{bwhi} is proved for $\mathcal{L}_{r_0}$ in a half-ball of unit size, scaling back we obtain \eqref{ineq_bwhi} for $\ld$, $u$, and balls with radius $R=r_0$, with a constant $C$ depending also on $r_0$. We can then cover $B_1^+$ with overlapping balls and semi-balls of size $r_0$, use that \eqref{ineq_bwhi} holds in each of these balls and a Harnack chain argument, to deduce \eqref{ineq_bwhi} in $B_1^+$ (see the proof of Theorem~2.1 in \cite{SS} for such a Harnack chain argument).

To prove Theorem \ref{bwhi}  for $\mathcal{L}_{r_0}$ in $B_1^+$ we use Lemma \ref{deggrowth}, which we proved above. Actually, once Lemma \ref{deggrowth} is available, Theorem \ref{bwhi} follows as in (the non-divergence case) \cite{Sir1}, the proof becomes essentially independent of the nature of the PDE. We sketch the argument, for completeness and convenience.

The inequality \eqref{ineq_bwhi} in Theorem \ref{bwhi} follows from Lemma \ref{deggrowth} in exactly the same way as \cite[Theorem 1.2]{Sir1} follows from \cite[Lemma 4.1]{Sir1}. We may repeat almost verbatim the argument on pages 7475-7478 in that paper. We note this argument uses only standard analysis, cube decomposition, the interior Harnack inequality and Lemma \ref{deggrowth}, so is independent of the nature (divergence or non-divergence) of the PDE.

To prove \eqref{ineq_bhi} in Theorem \ref{bwhi} we first observe that $u$ is H\"older continuous in $B_{3/2}^+$ (see Proposition \ref{holde} below). Further, we have the following Lipschitz estimate: if $u\in H^1(B_2^+)\cap C(\overline{B_2^+})$ is a weak solution of $\ldu\le f+\mathrm{div}(g)$ in $B_2^+$ with $u\le 0$ on $B_2^0$ then
\begin{equation}\label{lip}
u(x)\le C(\sup_{B_{3/2}^+} u + \|f\|_{L^q(B_2^+)}+  \|g\|_{L^\infty(B_2^+)}+\mathcal{M}_g(B_2^+))\,x_n,\quad x\in B_1^+.
\end{equation}
This can be proved in our setting in exactly the same way as \cite[Theorem 2.3]{Sir1}, by replacing the ABP inequality in the proof there by the Stampacchia generalized maximum principle (\cite[Theorem 8.15-8.16]{GT}) here, and by using the $C^1$-estimate in Theorem \ref{dek} together with the solvability results we already quoted (\cite[Theorems 3.1 and 3.3]{St}) in a sufficiently small neighborhood of the boundary.

Then, thanks to \eqref{lip} we can repeat the proof of
\cite[Theorem 1.3]{Sir1} (see page 7464 in \cite{Sir1}) in order to show that, in our divergence setting,
each weak solution of $\ldu\le f+\mathrm{div}(g)$ in $B_2^+$ with $u\le 0$ on $B_2^0$ is such that for each $r>0$
$$
\sup_{B_1} \left(\frac{u^+}{x_n}\right)\le C\left(\left(\int_{B_{3/2}} (u^+)^r\right)^{1/r} + \|f\|_{L^q(B_2)}+  \|g\|_{L^\infty(B_2)}+\mathcal{M}_g(B_2)\right).
$$
This, together with \eqref{ineq_bwhi}, gives \eqref{ineq_bhi}.

The third statement in Theorem \ref{bwhi} follows from \eqref{ineq_bwhi}, \eqref{ineq_bhi}, the interior (weak) Harnack inequallity,  and the same  covering/Harnack chain argument as above.

\subsection{Proof of Theorem \ref{kryl}}

Theorem \ref{kryl} follows from Theorem \ref{bwhi} and the method developed in \cite{SilS}. More precisely, we will see that Theorem \ref{kryl} is to the divergence case what \cite[Theorem 1.1]{SilS} is to the non-divergence one. For the reader's convenience we are going to fully review the proof and make explicit all parallels between the two works, giving all details at points where differences appear.

We recall we can assume that the boundary of $\Omega$ is locally flat, included in $\{x:x_n=0\}$. We first observe that global H\"older estimates are available for functions which satisfy the hypothesis of Theorem \ref{kryl} (this will replace the use of \cite[Proposition 2.6]{SilS} in our case).
\begin{prop}\label{holde} Assume $\ld^{(1)}, \ld^{(2)}$ are operators in the form \eqref{defdiv} which satisfy $(H1)$-$(H2)$. Assume $u\in H^1(B_{1}^+)$ is such that
 $$\ld^{(1)}[u]\le f^{(1)}+\mathrm{div}(g^{(1)}),\quad \ld^{(2)}[u]\ge f^{(2)}+\mathrm{div}(g^{(2)})\mbox{ in }B_{1}^+,\quad  u|_{B_{1}^0}\in C^{\overline \alpha}(B_{1}^0),$$
  for some $ \overline \alpha>0$. Then $u\in C^\alpha(B_{3/4}^+)$ and
\begin{equation}\label{holde1}
\|u\|_{C^\alpha(B_{3/4}^+)} \le C ( \|u\|_{L^\infty(B_1^+)} + \|f\|_{L^{p}(B_1^+)}+ \|g\|_{L^{q}(B_1^+)})
\end{equation}
for $\alpha, C>0$ depending on $n, \lambda,\Lambda, p,q,\overline \alpha$ and upper bounds on the Lebesgue norms of $\beta,b,c$ from $(H2)$ in $B_1^+$.
\end{prop}

\noindent{\it Proof.} This is \cite[Theorem 8.29]{GT} together with the remark at the end of Section 8.10 of \cite{GT} or \cite[Corollary 6.1]{Tr2}. Note these results were stated for solutions of equations, but their proofs were actually done for any function $u$ which is a supersolution and a subsolution of two different equations - since the proof is essentially an application of the WHI and bWHI to functions in the form $\sup_B u - u$ and $u-\inf_B u$ for well chosen balls and half-balls $B$. \hfill $\Box$
\smallskip

Therefore from now on we can assume without restricting the generality that $u\in H^1(B_{1}^+)\cap C^\alpha(B_{1}^+)$ as well as $c^{(i)}=\beta^{(i)}=0$ in Theorem \ref{kryl} (by replacing $f^{(i)}$ by $f^{(i)}-c^{(i)}u$, $g^{(i)}$ by $g^{(i)}+\beta^{(i)}u$).\smallskip

To parallel the notations in \cite{SilS} we set
$$
K=\sum_i
\|b^{(i)}\|_{L^q(B_{1}^+)} , \qquad L=\sum_i (
\|f^{(i)}\|_{L^q(B_{1}^+)}+\|g^{(i)}\|_{C^{\overline \alpha}(B_{1}^+)}).
$$

The inequality \eqref{ineq_bhi} plays a pivotal role in the proof of Theorem \ref{kryl}, providing the same bound as in \cite[Proposition 2.5]{SilS} (we will take $\varepsilon_0=1$ in that proposition). As before $e=(0,\ldots,0,1/2)$.
\begin{prop}\label{hopfe} Assume $u$ is as in Theorem  \ref{kryl}, with  $u\ge 0$ in $B_{1}^+$, $u=0$ on $B_0^1$, and $K\le 1$. Then
\begin{equation}\label{hopfe1}
u(x)\ge c_0(u(e)-C_0L)\, x_n, \quad x\in B_{3/4}^+,
\end{equation}
where  $c_0,C_0>0$ depend on $n, \lambda,\Lambda, q$, $\overline \alpha$, $\|A^{(i)}\|_{C^{\overline \alpha}(B_{1}^+)}$.
\end{prop}

\noindent{\it Proof.} This follows directly from \eqref{ineq_bhi}.\hfill $\Box$
\bigskip

Below we will also use  the Lipschitz estimate \eqref{lip} which replaces in our setting Proposition 2.4 in \cite{SilS}.
\begin{prop}\label{lipe} Assume $u$ is as in Theorem  \ref{kryl}, with   $u=0$ on $B_0^1$. Then
\begin{equation}\label{lipe1}
|u(x)|\le \bar C(\|u\|_{L^\infty(B_1^+)}+L)\, x_n, \quad x\in B_{3/4}^+,
\end{equation}
where  $\bar C>0$ depends on $n, \lambda,\Lambda, q$, $\overline \alpha$, $\|A^{(i)}\|_{C^{\overline \alpha}(B_{1}^+)}$, $K$.
\end{prop}
\noindent{\it Proof.} This follows directly from \eqref{lip} applied to $u$ and $-u$.\hfill $\Box$
\bigskip

\noindent{\bf Proof of Theorem \ref{kryl}}.  The first (and key) step in the proof is the particular case when $u$ vanishes on the boundary, $u|_{B_{1}^0}=0$.

\begin{prop}\label{vanish} Theorem \ref{kryl} is valid under the additional assumption\\ $u|_{B_{1}^0}=0$.
\end{prop}

\noindent {\it Proof.} Substituting $u$ by $u/W$, we can assume $W\le1$, that is, $|u|+L\le1$.

First, we observe that it is enough to establish the existence of a ``boundary gradient" at one point, for instance, it is enough to prove that there is a constant $G_0\in \mathbb{R}$, $|G_0|\le C$, such that \eqref{e:targetC1a} holds for $x_0=0$, $G(0)=(0,\ldots,0,G_0)$. Once this is proved, Proposition \ref{vanish} follows from a simple translation and analysis argument, described in the proof of \cite[Theorem 3.2]{SilS} (note that argument does not depend on the PDE).

So let us prove that such a constant $G_0$ exists. The proof uses the same idea as the one of \cite[Lemma 3.1]{SilS}, which is the corresponding result for non-divergence form equations (setting $A=G_0$ there). We can repeat the beginning of that proof, defining the sequences $r_k= 2^{-k}$, and $U_k$, $V_k$, to be built in such a way that $U_1 = \bar C$, $V_1= -\bar C$ ($\bar C\ge1$ is the constant from the Lipschitz estimate \eqref{lipe1}), $U_k$ is decreasing, $V_k$ is increasing,
$$
V_kx_n\le u(x)\le U_k x_n\;\mbox{ in } B_{r_k}, \qquad U_k-V_k \le Mr_k^{\alpha},
$$
which guarantees that the limit $G_0 = \lim_{k\to\infty} V_k=\lim_{k\to\infty} U_k$ is such that $|G_0|\le \bar C$ and has the desired property
$$
|u(x)-G_0x_n|\le 2^\alpha M|x|^{1+\alpha}, \qquad x\in B_{1/2}^+.
$$
The constants $\alpha, M$ are defined by
$$\alpha=\frac{1}{2}\min\{ 1-n/q, \bar \alpha, \tilde\alpha\},\quad\mbox{ with }\; 2^{-\tilde\alpha}= 1-c_0/8,
$$
where $c_0\le1$ is the constant from \eqref{hopfe1} and
$$
M=2\bar Cr_{k_0}^{-\alpha}=\bar C 2^{1+\alpha k_0} ,\qquad \mbox{with }\; k_0 \mbox{ chosen so that }
$$
$$
(\bar CK+1)r_{k_0}^{1-\alpha-n/q}< 1/(16C_0), \qquad (\bar C\|A^{(i)}\|_{C^{\overline \alpha}(B_{1}^+)}+1)r_{k_0}^{\bar\alpha-\alpha}< 1/(16C_0),
$$
 where $C_0\ge1$ is the constant from \eqref{hopfe1}.

Exactly as in \cite[Lemma 3.1]{SilS} we define $U_1=\ldots=U_{k_0}$, $V_1=\ldots=V_{k_0}$ and construct iteratively $U_k,V_k$ for $k>k_0$. The iterative procedure is the same,  up to the introduction of the rescaled function
\begin{equation}\label{noy1}
v_k(x) = r_k^{-1-\alpha}(u(r_kx) - V_kr_kx_n), \quad x\in B_1^+.
\end{equation}
where the only real difference with the proof of \cite[Lemma 3.1]{SilS} appears: in our setting the new function $v_k$ satisfies in $B_1^+$
\begin{align*}
-\frac{1}{r_k}\mathrm{div}(A(r_kx)(r_k^\alpha Dv_k(x) + V_k e_n)) + b(r_kx)(r_k^\alpha Dv_k(x) + V_k e_n))\\ \le(\ge) f(r_k x) + \frac{1}{r_k}\mathrm{div}(g(r_kx))
\end{align*}
(we do not write the indices $i=1,2$ for display convenience), which can be rewritten as
\begin{align*}
-\mathrm{div}(A(r_kx)Dv_k(x)) + r_kb(r_kx)Dv_k(x)
\le(\ge) r_k^{1-\alpha}( V_k b(r_kx) e_n +  f(r_k x) )\\ + \mathrm{div}\left(V_k\frac{A(r_kx)-A(0)}{r_k^\alpha}e_n\right) + \mathrm{div}\left(\frac{g(r_kx)-g(0)}{r_k^\alpha}\right).
\end{align*}
We observe that
$$
\|A(r_kx)\|_{C^{\bar\alpha}(B_1^+)}\le \|A(x)\|_{C^{\bar\alpha}(B_1^+)},
$$
$$r_k\|  b(r_kx)\|_{L^q(B_1^+)}\le  r_k^{1-n/q}\|b\|_{L^q(B_1^+)} \le K r_k^{1-n/q},
$$
$$
r_k^{1-\alpha}\| V_k b(r_kx)\|_{L^q(B_1^+)}\le \overline{C} r_k^{1-\alpha-n/q}\|b\|_{L^q(B_1^+)}\le \overline{C} K r_k^{1-\alpha-n/q},
$$
$$
r_k^{1-\alpha}\| f(r_kx)\|_{L^q(B_1^+)}\le r_k^{1-\alpha-n/q}\|f\|_{L^q(B_1^+)}\le   r_k^{1-\alpha-n/q}, \quad (L\le 1),
$$
$$
\left\|V_k\frac{A(r_kx)-A(0)}{r_k^\alpha}\right\|_{C^{\bar\alpha}(B_1^+)}\le  \overline{C}\|A\|_{C^{\bar\alpha}(B_1^+)} r_k^{\bar\alpha-\alpha}, $$
$$\left\|\frac{g(r_kx)-g(0)}{r_k^\alpha}\right\|_{C^{\bar\alpha}(B_1^+)}\le  \|g\|_{C^{\bar\alpha}(B_1^+)}r_k^{\bar\alpha-\alpha}\le r_k^{\bar\alpha-\alpha},\quad (L\le 1),
$$

Hence  the last differential inequality is in the form
$$
-\mathrm{div}(A(r_kx)Dv_k(x)) + \hat b_k  Dv_k \le(\ge) \hat f_k + \mathrm{div}(\hat g_k)
$$
where $\hat f_k$, $\hat g_k$ are such that the corresponding quantities $\hat K$, $\hat L$ satisfy $\hat K< 1$ and $\hat L< 1/(8C_0)$, thanks to the choice of $k_0$ and $k\ge k_0$. So we can apply Proposition \ref{hopfe}, getting
$$
v_k(x)\ge c_0(v_k(e) - C_0\hat L)x_n, \qquad x\in B_{1/2}^+.
$$
This inequality replaces the inequality (3.4) in \cite{SilS} and the rest of the proof of Lemma 3.1 there is repeated without any changes\footnote{We note a small misprint in \cite{SilS}, one sets there $\varepsilon_1=1/(8C_0)$.}, the PDE is no longer used.

 Thus  Proposition \ref{vanish} is established.\hfill $\Box$
\medskip

To extend Proposition \ref{vanish} to arbitrary boundary data, we can just remove a $C^{1,\bar\alpha}$-function from $u$. For instance, for $x=(x^\prime,x_n)\in B_1$ we set $$\psi(x) = u(x^\prime,0) \qquad\mbox{and}\qquad v= u-\psi.$$ Then $v$ vanishes on the flat boundary and solves the same equation as $u$, with $f$ replaced by $f-bD\psi$ and $g$ replaced by $g+AD\psi$. Since obviously $\|\psi\|_{C^{1,\bar\alpha}(B_1^+)} = \|u\|_{C^{1,\bar\alpha}(B_1^0)}$, it is trivial to see that Proposition  \ref{vanish} applied to~$v$ implies the statement of Theorem \ref{kryl}.\hfill $\Box$

\bigskip

\noindent{\bf Funding acknowledgements}. B. Sirakov is supported by CNPq grant 310989/2018-3 and FAPERJ grant E-26/201.205/2021. M. Soares is supported by a post-doctoral grant DGAPA--UNAM.\bigskip



{\upshape
Fiorella Rend\'on\\
Departamento de Matem\'atica,\\
Pontif\'{\i}cia Universidade Cat\'olica do Rio de Janeiro (PUC-Rio)\\
22451-900, G\'avea, Rio de Janeiro-RJ, Brazil\\
fiorellareg@gmail.com\bigskip

Boyan Sirakov\\
Departamento de Matem\'atica,\\
Pontif\'{\i}cia Universidade Cat\'olica do Rio de Janeiro (PUC-Rio)\\
22451-900, G\'avea, Rio de Janeiro-RJ, Brazil\\
bsirakov@puc-rio.br\bigskip

Mayra Soares\\
Instituto de Matem\'atica e Estat\'{\i}stica\\
Universidade Federal de Goi\'as\\
Campus Samambaia, Avenida Esperança, s/n.\\
Cep: 74690-900, Goi\^ania, Goi\'as, Brazil\\
mayra.soares@ufg.br
}

\end{document}